\documentclass[11pt]{article}

\usepackage{amsfonts,amsmath,amsthm,amstext}

\newenvironment{map_machine}[5]
   {
    
    \begin{array}[#1]{rcl}
   #2 & {-}\!\!{-}\!\!{\longrightarrow} & #3\\ 
   #4 & \mapstochar\kern-.25em{-}\!\!{\longrightarrow} & #5
   }{\end{array}}

\newcommand{\compress} 
  {\setlength{\partopsep}{0in}
  \setlength{\itemsep}{0in} \setlength{\topsep}{0in}
  \setlength{\parsep}{0in} \setlength{\parskip}{0in}}
  
\newcommand{\romanate}
{\renewcommand{\labelenumi}{{\rm(\roman{enumi})}}
 \renewcommand{\labelenumii}{{\rm(\alph{enumii})}}}

\def\implies{\ensuremath{\Rightarrow}}

\newcommand{\define}[1]{\textbf{#1}}

\renewcommand{\phi}{\varphi}
\renewcommand{\tilde}{\widetilde}
\renewcommand{\epsilon}{\varepsilon}

\newcommand{\C}[1]{C_{#1}^{}}
\newcommand{\genby}[1]{\langle#1\rangle}
\newcommand{\rtsys}{\Phi}
\newcommand{\rtbase}{\Delta}
\newcommand{\st}{\,|\,}
\newcommand{\inv}{{}^{-1}}

\newcommand{\two}{{}^{(2)}}
\newcommand{\nottwo}{{}^{(\ne2)}}

\newcommand{\rto}[1]{\stackrel{\displaystyle #1}{{-}\kern-.25em{-}\!\!{\longrightarrow}}}
\newcommand{\restrict}[1]{\big|_{\scriptstyle #1}^{}}

\newcommand{\Z}{\mathbb{Z}}
\newcommand{\CC}{\mathbb{C}}

\DeclareMathOperator{\BCpairs}{BC-pairs}
\DeclareMathOperator{\height}{ht}
\DeclareMathOperator{\rank}{rank}
\DeclareMathOperator{\orthog}{O}
\renewcommand{\O}{\orthog}

\DeclareMathOperator{\GL}{GL}
\DeclareMathOperator{\Cl}{Cl}
\DeclareMathOperator{\symplectic}{Sp}
\def\Sp{\symplectic}
\DeclareMathOperator{\SO}{SO}

\newcommand{\deven}[1]%
  {\ensuremath{\left|\left\{\begin{array}{@{}l@{}}\text{even, %
  singular}\\[-.5ex]%
  \text{parts of }\lambda\end{array}\right\}\right|}}
   
\newtheorem{lemma}{Lemma}[section]              
\newtheorem{cor}[lemma]{Corollary}
\newtheorem{theorem}{Theorem}

\newtheorem*{richardsons_theorem}{Richardson's Theorem}
\theoremstyle{definition}
\newtheorem{innerproof}[lemma]{Proof}
\newenvironment{numproof}{\begin{innerproof}}{\qed\end{innerproof}}

\newtheorem{InnerRemarks}[lemma]{Remarks}
\newcounter{RemarkCounter}
\newcommand{\rem}{\stepcounter{RemarkCounter}{\rm 
    (\alph{RemarkCounter})}\ }
\newenvironment{remarks}{\setcounter{RemarkCounter}{0}\begin{InnerRemarks}}{\end{InnerRemarks}}
\newtheorem{remark}[lemma]{Remark}

\newcounter{arrowlength}
\newcommand{\To}[1]{%
    \loop \relbar\joinrel\stepcounter{arrowlength}
    \ifnum\value{arrowlength}<#1\repeat
    \to\setcounter{arrowlength}{0}}

\newbox{\myaddress}
\savebox{\myaddress}{{\rm \normalsize 334 Hill Center, Rutgers 
University, Piscataway NJ, 08854}}
\newbox{\myemail}
\savebox{\myemail}{\tt \normalsize duck@math.rutgers.edu}

\title{Jordan Blocks of Richardson Classes in the Classical Groups and 
the Bala-Carter Theorem}

\author{W. Ethan Duckworth\\
\usebox{\myaddress}\\
\usebox{\myemail}}
\date{January 10, 2003}

\begin{document}
\maketitle

\begin{abstract}This paper provides new, relatively simple proofs of some
important results about unipotent classes in simple linear algebraic
groups.  We derive the formula for the Jordan blocks of the Richardson
class of a parabolic subgroup of a classical group.  This result was
originally due to Spaltenstein.  Secondly, we derive, for good
characteristic, the description of the natural partial order of
unipotent classes of a classical group in terms of their Jordan
blocks.  This result was originally due to Gerstenhaber and Hesselink. 
As a consequence we obtain a proof of the Bala-Carter Theorem which
holds even in certain bad characteristics (this proof requires the
prior classification of unipotent classes, unlike the original proofs
due to Bala, Carter and Pommerening).
\bigskip

\end{abstract}

\section{Introduction}
Let $G$ be a connected reductive linear algebraic group.  Richardson
\cite{richardson} made a vital contribution to the study of unipotent
classes in algebraic groups by associating to each parabolic subgroup
of $G$ a unipotent class of $G$.  This result has had surprisingly
powerful implications, some of which we will discuss below.  The
following theorem is one version of Richardson's result (see also
\cite{carter} or \cite{steinberg} for other proofs).

\begin{richardsons_theorem}
\label{Richardsons_theorem}
\index{Richardson's theorem}
Let $G$ be a connected reductive group,
$P$ a para\-bolic subgroup with unipotent radical $Q$ and Levi factor
$L$.  The following hold:
\begin{enumerate}
\romanate\compress

\item 
There exists a unique unipotent $G$-class $C$ such that $C\cap Q$ is
open and dense in $Q$.

\item 
$C\cap Q$ forms a single $P$-class.

\item
If $u\in C\cap Q$ then $\C G(u)^\circ = \C P(u)^\circ$, whence these
centralizers have dimension $\dim L$.

\item 
Let $Z$ be the center of $G$, let $Q'$ be the derived subgroup of $Q$. 
Then $\dim L/Z \ge \dim Q/Q'$.
\end{enumerate}
\end{richardsons_theorem}

Spaltenstein \cite{spaltenstein} studies generalizations of this
result to the case where $G$ is non-connected.  We will always have
$G$ connected except when $G=\O_n$.

We call $C$ the \define{Richardson class}\index{Richardson class} of
$P$ and we call $C\cap Q$ the \define{Richardson orbit}%
\index{Richardson orbit} in $Q$.  \bigskip

For many questions it is of fundamental importance to be able to find
the Jordan blocks of a unipotent class.  The next result indicates how 
to do this for Richardson classes, but first we introduce some 
standard notation.  

A \define{partition}\index{Partition} of $n$ is a sequence of natural
numbers which add to $n$; we assume that the sequence is weakly
decreasing unless indicated otherwise.  We write a partition
$\lambda$ as $(\lambda_1, \lambda_2, \dots)$ or $(\lambda_1\ge
\lambda_2 \ge \dots)$ or $(1^{c(1)}, 2^{c(2)}, 3^{c(3)},\dots)$ where
$c(x)$ is the multiplicity of $x$ in $\lambda$.  We call $\lambda_i$ a
\define{part} of $\lambda$.  
The \define{dual}\index{dual of a partition}\index{Partition!dual of}
of $\lambda$ is a partition of $n$ which we write as $\lambda^*$ and
which has parts defined as follows: $\lambda_i^*$ equals the number of
parts of $\lambda$ which are greater than or equal to $i$ (i.e.
$\lambda_i^*$ equals the number of indices $j$ such that $\lambda_j
\ge i$).


Let $G\in \{\SO_{2n}, \SO_{2n+1}, \Sp_{2n}\}$.  We fix a root base
$\rtbase$ for $G$ and label the nodes $\{\alpha_1,\dots,\alpha_n\}$ as
in \cite{bourbaki}.  Given a parabolic subgroup $P$ let $J\subseteq
\rtbase$ such that $P$ is conjugate to the standard parabolic
associated with $J$.  If $\alpha_n\in J$ then let $m$ be the largest
integer such that such that the last $m$ nodes $\alpha_{n-m+1},\dots,
\alpha_n$ are contained in $J$.  Then the Levi factor $L$ of $P$ can
be written as $L=\GL_{n_1}\dots\GL_{n_s} \Cl_m$ where $\Cl_m\in
\{\SO_{2m}, \SO_{2m+1}, \Sp_{2m}\}$.  We extend this notation also to
the cases $\alpha_n\not\in J$ and to $G=\GL_n$ by taking $\Cl_m=1$ and
$m=0$.  In this way each parabolic subgroup of $\GL_n$, $\SO_{2n}$,
$\SO_{2n+1}$, $\Sp_{2n}$ determines a partition $n= n_1+\dots +n_s
+m$.  We call this the \define{Levi partition}\index{Levi partition}
and write it either as $\Lambda = (n_1,\dots,n_s) \oplus m$ or
$\Lambda = (1^{c(1)}, 2^{c(2)}, \dots)\oplus m$ where $c(x)$ is the
multiplicity of the part $x$ in the $n_i$ and the notation
``$\oplus$'' indicates that we have an ordered pair consisting of the
partition $(n_1,\dots,n_s)$ and the number $m$.

\begin{table*}%
\caption{Jordan blocks of a Richardson class}%
\label{table_of_Jordan_blocks_of_Rich_class}%
$$\begin{array}{|rcll|}
\hline
&& \multicolumn{2}{l|}{G=\GL_n,\ \psi \text{ is the identity}}\\
\hline
%
%
&& \multicolumn{2}{l|}{G=\SO_{2n},\ p\ne 2}\\
\psi(m) & = & 2m &\\
\psi(j^{c(j)}) & = & j^{2c(j)} 
      & \text{if }j\text{ is even or } j\le 2m \\
\psi(j^{c(j)}) & = & j+1,\ j^{2c(j)-2},\ j-1  
      & \text{if } j\text{ is odd and }j> 2m\\
\hline
%
&& \multicolumn{2}{l|}{G = \SO_{2n},\ p=2}\\
\psi(m) & = & 2m & \\
\psi(j^{c(j)}) & = & j^{2c(j)}             
      & \text{if } j \text{ is even}\\
\psi(j^{c(j)}) & = & j+1,\ j^{2c(j)-2},\ j-1 
      & \text{if } j\text{ is odd and }j\le 2m\\
\psi(j^{c(j)}) & = & j+1,\ j-1              
      & \text{if }j \text{ is odd, } j> 2m\text{ and } c(j) = 1\\
\psi(j^{c(j)}) & = & (j+1)^2,\ j^{2c(j)-4},\ (j-1)^2  
      & \text{if }j\text{ is odd, }j> 2m\text{ and } c(j) \ge 2\\
\hline
&&\multicolumn{2}{l|}{G=\SO_{2n+1},\ p\ne2}\\
\psi(m) & = & 2m+1 &\\
\psi(j^{c(j)}) &=& j^{2c(j)}             
      & \text{if } j \text{ is odd and }j>2m+1 \\
\psi(j^{c(j)}) &=& j^{2c(j)}             
      & \text{if }j\le 2m+1\\
\psi(j^{c(j)}) &=& j+1,\ j^{2c(j)-2},\ j-1 
      & \text{if } j \text{ is even and }j>2m+1\\
\hline
&&\multicolumn{2}{l|}{G=\Sp_{2n}}\\
\psi(m) & = & 2m & \\
\psi(j^{c(j)}) &=& j^{2c(j)}             
      & \text{if } j \text{ is even or }j\ge 2m\\
\psi(j^{c(j)}) &=& j+1,\ j^{2c(j)-2},\ j-1 
      & \text{if }j \text{ is odd and }j<2m\\
\hline
\end{array}$$
\end{table*}
%
%
%

\begin{theorem}[{\cite[II.7.4]{spaltenstein}}]
\label{theorem_Jordan_blocks_of_rich_classes} 
Let $G$ be one of $\GL_n$, $\SO_{2n}$, $\SO_{2n+1}$, $\Sp_{2n}$ and
exclude the case $\SO_{2n+1}$ if $p=2$.  Let $P$ be a parabolic
subgroup of $G$, let $\Lambda$ be the Levi partition of $P$, define
the map $\psi$ as in
table~{\rm\ref{table_of_Jordan_blocks_of_Rich_class}}, and let
$\lambda$ be the partition of Jordan block sizes of the Richardson
class of $P$.  Then $\lambda$ equals $\psi(\Lambda)^*$, the dual of
$\psi(\Lambda)$.
\end{theorem}

\begin{remarks} 
\rem It is easy to extend this result to the case of
$G=\SO_{2n+1}$ and $p=2$.  One applies the formula for $\Sp_{2n}$
using the same Levi partition and adds one block of size $1$ to the
result.  
\rem Spaltenstein's formulas appear rather different from those
presented here (and have some minor mistakes).  (Spaltenstein also
determines the index $\epsilon$ in his notation or, equivalently, the
singularity of the parts of $\lambda$.  See \cite{duckworth} for a
discussion of this notation.)
%
\end{remarks}

Let $C_1$ and $C_2$ be two unipotent classes of $G$.  We
define $C_1 \le C_2$ if and only if $C_1\subseteq \overline{C_2}$
(where $\overline {C_2}$ is the closure of $C_2$).   This is the 
natural partial order on unipotent classes.  

Let $\lambda = (\lambda_1, \lambda_2,\dots)$ and
$\mu=(\mu_1,\mu_2,\dots)$ be two partitions.  We define $\lambda
\le \mu$ if and only if for all $j\ge 1$ we have $\sum_{i=1}^j
\lambda_i \le \sum_{i=1}^j \mu_i$.  This is the \define{dominance} 
partial order on partitions.  

\begin{theorem}[\cite{hesselink2}, \cite{gerstenhaber}, 
{\cite[I.2.10]{spaltenstein}}]
\label{theorem_on_order_relations}
Let $G\in\{\GL_n, \O_n, \SO_n,\Sp_n\}$.  Let $C_\lambda$ and $C_\mu$
be two unipotent classes in $G$ with $\lambda$ and $\mu$ the
partitions consisting of the Jordan blocks of $C_\lambda$ and $C_\mu$
respectively.  Assume either that $p\ne2$ if $G\ne \GL_n$ or that
$\mu$ has no even parts with even multiplicity. Then $\lambda< \mu$
if and only if $C_\lambda < C_\mu$.
\end{theorem}

Spaltenstein \cite[I.2.10]{spaltenstein} generalizes this result to
all unipotent classes in bad characteristics, but we will not discuss
his generalization here.
%
%
%
%
%
%

We now introduce the necessary terminology to state the Bala-Carter
Theorem which gives a parameterization of unipotent classes in
simple algebraic groups. 

Let $G$ be a connected reductive algebraic group with root system
$\rtsys$ and root base $\rtbase$.  Fix $J\subseteq \rtbase$ and let
$P$ be the standard parabolic subgroup corresponding to $J$.  Let
$\beta\in \rtsys$ and write $\beta = \sum_{\alpha\in
\rtbase}n_{\alpha}\, \alpha$.  The $P$-height is defined to be
$\height_P^{}(\beta) = \sum_{\alpha\in \rtbase -J}^{} n_{\alpha}$.

Let $L$ be a Levi factor for $P$, $Q$ the unipotent radical of $P$,
and $\rtsys(Q)$ the roots of $Q$.  We say $P$ is
\define{distinguished}\index{distinguished parabolic}\index{Parabolic
subgroup!distinguished} if $\dim L/Z(G)$ equals the number of roots in
$\rtsys(Q)$ with $P$-height equal to $1$.

If $Q'$ is the derived subgroup of $Q$ then Richardson's Theorem (iv)
implies that $\dim L/Z(G)\ge \dim Q/Q'$ for all $P$.  If $P$ is
distinguished then
$$\dim L/Z(G) = \dim Q/Q'\makebox[0in][l]{\hspace*{1in} $(*)$}.$$  
The converse holds provided $p\ne 2$ if the Dynkin diagram of $G$
contains double bonds, and $p\ne 3$ if the Dynkin diagram of $G$
contains triple bonds (see \cite{bala-carter} or \cite{ABS}).

The work in \cite{bala-carter}, \cite{carter} etc.\ takes condition
$(*)$ as the definition of distinguished, but then applies this
definition only with the restrictions on $p$ just described.
Thus, the definition we have given here takes the usual list of 
distinguished parabolics and uses this same list even when $p$ 
equals $2$ or $3$.  We refer the reader to \cite{carter} for a 
list of the distinguished parabolics (note however that there is a 
mistake in the second formula for $D_n$).  
%

Throughout this paper, a \define{Levi subgroup} means a Levi factor of
a parabolic subgroup.  Let $L$ be a Levi subgroup of $G$ and $u\in L$
a unipotent element.  Then $u$ is \define{distinguished} in $L$ if $u$
is not contained in any proper Levi subgroup of $L$.  If $L=G$ has
trivial center this is equivalent to having $\C G(u)$ contain no
non-trivial torus (see Lemma~\ref{min_Levi_unique_up_to_conj} below).

For a reductive group $G$ let $\BCpairs(G)$ denote the pairs $(L,P)$
where $L$ is a Levi subgroup and $P$ is a distinguished parabolic
subgroup of $L$.  Let $\psi$ (or $\psi_G^{}$) denote the map from such
pairs to unipotent classes in $G$ obtained by extending the Richardson
class of $P$ (in $L$) to a $G$-class.

\begin{theorem}[Bala-Carter \cite{bala-carter}, Pommerening 
\cite{pommerening}]
\label{bala-carter-theorem}\index{Bala-Carter theorem}
Let $G$ be a simple algebraic group and let $\psi= \psi_G^{}$ be as just defined.
The following hold:
\begin{enumerate}
\romanate\compress

\item 
If $X$ is a Levi subgroup the following diagram commutes: 
$$\begin{array}{ccc} 
X\text{-classes in }\BCpairs(X) & \rto{\psi_X^{}} & \text{unipotent classes in }X\\
\downarrow & \circ & \downarrow\\
G\text{-classes in }\BCpairs(G) & \rto{\psi} & \text{uniptotent classes in }G
\end{array}$$
where the vertical maps extend an $X$-class to the corresponding
$G$-class.

\item 
Let  $\psi(L,P) = C$ and $u\in C\cap L$.  Then $u$ is distinguished 
in $L$.

\item 
The map $\psi$ is injective.  It is a bijection except in the
following cases: $G\in \{B_n,C_n,D_n\}$ and $p=2$; $(G,p)$ is one 
of $(E_7,2)$, $(E_8,2)$, $(E_8,3)$, $(F_4,2)$ or $(G_2,3)$ in which 
cases there are $1$, $4$, $1$, $4$ and $1$ extra classes 
respectively.  
\end{enumerate}
\end{theorem}

\begin{remarks}
\label{remarks_after_bala_carter}
Although part (i) is obvious, we state it here to bring 
attention to some of the following points. 
\rem Part (i) makes the Bala-Carter Theorem more useful than Jordan
blocks for comparing unipotent classes in $X$ and unipotent classes in
$G$.  For example, let $G=E_6$, $X=D_5T_1$ and let $C$ be the
unipotent class of $X$ which has two Jordan blocks of size $5$ in the
natural module for $\SO_{10}$.  The Jordan blocks do not make it clear
which class $C$ corresponds to in $G$.  However, the Bala-Carter label
for $C$ is $A_4$ (i.e. $C$ is represented by a regular element of a
Levi subgroup of type $A_{4}$) both as a class of $X$ and when it is
extended to a class of $G$.
\rem If $X$ is a maximal rank reductive subgroup which is not a Levi
subgroup, one may often still obtain a commutative diagram similar to
that in part (i).  For instance let $G=E_6$, $X=A_2A_2A_2$ and
$(L,P)\in \BCpairs(X)$ where $L$ is a proper Levi subgroup of $X$. 
Then $(L,P)\in \BCpairs(G)$ and the same result is obtained if one
first extends $(L,P)$ to a $G$-class and then takes the unipotent
$G$-class, or if one first takes the unipotent $X$-class and then
extends this to a $G$-class.
\rem Parts (i) and (iii) show that in most cases the intersection of a
unipotent $G$-class with a Levi subgroup forms a single unipotent
class for the Levi subgroup.  If this is not the case then $G=E_r$,
$L$ is of type $D_n$ and the unipotent class is of type $A_{n-1}$.  
\rem Part (iii) is a stronger version of the Bala-Carter-Pommerening
Theorem than usually appears (as in the references above or
\cite{carter}, \cite{humphreys3}), although this version seems to be
known or assumed by specialists in the field (see, for example,
\cite{lawther2}).  In addition, the proof given in this paper (see
Proof~\ref{proof_of_bala_carter}) uses the classification of unipotent
classes for each simple algebraic group whereas the standard proof (as
in \cite{carter}) constructs the inverse of $\psi$ (at the level of
the Lie algebra) and is independent of these classifications.
\end{remarks}

\section{Recollections and Conventions}
\label{conventions}
All algebraic groups in this paper are affine and defined over a fixed
algebraically closed field of characteristic $p\ge 0$.  

The groups $\SO_{2n}$, $\SO_{2n+1}$ and $\Sp_{2n}$ are defined in
terms of a bilinear form and a quadratic form which we will usually
denote by $\beta$ and $\phi$ respectively.  Let $V$ be the natural
module for one of these groups.  A subspace $W$ is \define{totally
singular} if $\phi\restrict W$ is identically zero (which implies that
$\beta\restrict {W\times W}$ also equals zero); 
it is \define{nonsingular} if $\beta\restrict{W}$ has trivial radical. 
If $G=\GL_n$ we consider each subspace of its natural module to be
totally singular.  If $W$ is a nonsingular subspace then $\Cl(W)$
denotes the classical group of the same type as $G$ defined on $W$.
%

Let $G$ be a classical group with natural module $V$.  A \define{flag}
is a sequence of nested subspaces.  Let $f$ be the flag $ W_0 \le W_1
\le \dots \le W_\ell = V$.  Then $f$ has length $\ell$ and is
\define{totally singular} if for each $i$ either $W_i$ is totally
singular or $W_i=W^\perp$ for some totally singular subspace $W \le V$
(if $V$ is nonsingular this is equivalent to requiring that either
$W_i$ or $W_i^\perp$ be totally singular).  A subgroup of $G$ is
parabolic if and only if it is the stabilizer of a totally singular
flag.  Let $L = \GL_{n_1} \dots \GL_{n_s} \Cl_m$ be the Levi factor of
a parabolic subgroup $P$.  We say a flag $f$ of length $\ell$ is a
\define{natural flag} for $P$ if the following hold: $f$ is totally
singular, $P$ is the stabilizer of $f$, $\ell = s$ if $G=\GL_n$, $\ell
= 2s$ if $m=0$ and $G\in \{\SO_{2n}, \Sp_{2n}\}$, and $\ell = 2s+1$ if
$m\ge 1$ or $G = \SO_{2n+1}$.  The unipotent radical of a parabolic 
equals the set of elements which act trivially upon each factor in a 
natural flag.

In the classical groups the unipotent classes are described using
partitions.  We will mention only a few facts here and refer the
reader to \cite{carter} or \cite{duckworth} for more complete
information.  Let $G$ be one of $\GL_n$, $\O_n$, $\SO_{n}$, $\Sp_{n}$,
let $C$ be a unipotent class of $G$ and let $\lambda$ be the partition
of $n$ consisting of the Jordan block sizes of $C$.   

The \define{parity conditions} on $\lambda$ refer to the following
requirements: if $G \in \{\O_n, \SO_n\}$ and $p\ne 2$ then each even
part of $\lambda$ must have even multiplicity; if $G=\Sp_{n}$ or
$G\in\{\O_n,\SO_n\}$ and $p=2$ then each odd part of $\lambda$ must
have even multiplicity; if $G=\SO_n$ with $n$ even then $\lambda$ must
have an even number of parts.

If $G\in \{\O_n, \Sp_n\}$ and $\lambda$ has no even parts with even
multiplicity then all unipotent elements with Jordan blocks equal to
$\lambda$ form a single $G$-class.  This is generally not the case if 
$G\ne \GL_n$ and $p=2$.  

If $u\in C$ we say a part $x$ of $\lambda$ is nonsingular if there
exists a Jordan chain of $u$ (i.e. a sequence of vectors
$(v_i)_{i=0}^x$ such that $v_0=0$ and $(u-1)v_i = v_{i-1}$ for all
$i\ge 1$) which generates an $x$-dimensional nonsingular subspace.

\begin{lemma}[\cite{duckworth}]
\label{jordan_chain_degeneracy} Let $G$ be $\O_n$, $\SO_n$ or $\Sp_n$,
with natural module $V$, bilinear form $\beta$, $u\in G$ a unipotent
element, $\lambda$ the Jordan blocks of $u$, $x$ a part of $\lambda$
and $v_1, \dots, v_x$ a Jordan chain of $u$.
\begin{enumerate}
\romanate\compress

\item The subspace $\genby{v_1,\dots,v_x}$ is nonsingular if and only
if $\beta(v_i,v_j)\ne 0$ for some or, equivalently, for all $\,i,j>0$
with $i + j = x + 1$.  If $V$ is nonsingular and $x$ has multiplicity
$1$ then the subspace $\genby{v_1,\dots,v_x}$ is nonsingular.

\item If $G$ equals $\O_n$ or $\SO_n$ and $p\ne 2$ then $x$ is
nonsingular if and only if $x$ is odd.  If $G=\Sp_n$ and $p\ne 2$ then
$x$ is nonsingular if and only if $x$ is even.  In any case, if $x\ne
1$ and the multiplicity of $x$ is odd then $x$ is nonsingular.
%
%
%
\end{enumerate}
\end{lemma}

\begin{remarks}
\label{remark_on_spaltenstein_centralizer_formula}
\rem
A stronger statement than given here is possible.  In particular, 
keeping track of information about singularity of partitions is enough to 
parameterize the unipotent classes of $\O_n$ and $\Sp_n$ in 
characteristic $2$.  
\rem Spaltenstein \cite[I.2.8]{spaltenstein} gives the following
expression (using different terminology) for the dimension of the
centralizer of a unipotent element.  Suppose $p=2$ and $G$ equals $\Sp_n$ or
$\SO_n$ with $n$ even..  Let $u\in G$ be a unipotent element
and let $u_\CC^{} \in \Sp_n(\CC)$ be a unipotent element with the same
Jordan blocks as $u$.  Then $\dim \C G(u)$ equals $\dim \C
{\Sp_n(\CC)}(u_\CC^{})$ plus the number of even, singular parts of
$\lambda$.
\end{remarks}

If $G$ is a reductive (not necessarily connected) group, a unipotent
element is \define{regular}\index{Unipotent class!regular elements} if
the dimension of its centralizer equals the rank of $G$.  If $G$ is
connected then the regular elements form a single unipotent class (see
\cite{steinberg_regular_elements} or \cite{carter}), which is the
Richardson class of the Borel subgroups.  If $G$ is not connected then
the number of regular classes is at most the number of connected
components (see \cite{spaltenstein} for more on this and the
connection with Richardson classes).

\begin{lemma}
\label{lemma_for_regular_elements}
\index{Unipotent class!Jordan blocks of regular elements}
\index{Regular unipotent elements} Let $G$ be one of $\GL_n$,
$\SO_{2n+1}$, $\SO_{2n}$, $\Sp_{2n}$ or $\O_{2n}$ and exclude the case
$\SO_{2n+1}$ with $p=2$.  Let $\lambda$ be the Jordan blocks of a
regular unipotent class.  If $G=\GL_n$ then $\lambda = n$.  If
$G=\SO_{2n+1}$ then $\lambda = 2n+1$.  If $G=\SO_{2n}$ then $\lambda$
equals $(2n-1,1)$ or $(2n-2,2)$ according as $p\ne2$ or $p=2$
respectively.  If $G=\Sp_{2n}$ then $\lambda = 2n$.  If $G=\O_{2n}$
and $p=2$ then there are two regular unipotent classes and these have
Jordan blocks of sizes $2n$ and $(2n-2,2)$.  In all cases all parts of
$\lambda$ are nonsingular.
\end{lemma}

\begin{proof}
This follows from an easy dimension calculation.
\end{proof}

\section{Proof of Theorems
\ref{theorem_Jordan_blocks_of_rich_classes} and
\ref{theorem_on_order_relations}}

%
%
%
%

For this section we use the following notation and assumptions (with
three explicit exceptions marked by the phrase ``Contrary to our usual
assumptions \ldots'').  We assume throughout that $G$ is one of
$\GL_n$, $\SO_{2n+1}$, $\SO_{2n}$, $\Sp_{2n}$ and exclude the case
$\SO_{2n+1}$ when $p=2$.  Let $V$ be the natural module for $G$ and
$\beta$ the bilinear form on $V$ if $G\ne \GL_n$.  Let $P$ be a proper
parabolic subgroup (we allow $P=G$ in the statement of
Theorem~\ref{theorem_Jordan_blocks_of_rich_classes}, but if this holds
there is nothing to prove).  Let $Q$ be the unipotent radical of $P$
and $f = (0 = W_0 < \dots < W_\ell = V)$ a natural flag.  Let $\Lambda
= (n_1,\dots,n_s)\oplus m = (1^{c(1)}, 2^{c(2)},\dots)\oplus m$ be the
Levi partition of $P$.

For any $g\in Q$ we let $\lambda(g) = (\lambda_1(g),
\lambda_2(g),\dots)$ be the partition of Jordan blocks of $g$.  We fix
$u\in Q$ which represents the Richardson orbit in $Q$.  We fix
$\lambda = (\lambda_1, \lambda_2,\dots) = \lambda(u)$ and $\mu =
\psi(\Lambda)^*$.  We wish to prove that $\lambda = \mu$.

Essentially the proof of
Theorem~\ref{theorem_Jordan_blocks_of_rich_classes} is inductive.  We
will produce the largest one or two Jordan blocks of $\lambda$ and
show that they equal the largest one or two parts of $\mu$ and that
they generate a nonsingular subspace $V_r$.  We will then look at the
action of $u$ on $V/V_r$ and induct.  See
Remark~\ref{remark_gather_pieces_for_main_proof} for the main steps in
this proof.

Theorem~\ref{theorem_on_order_relations} is proven in
Corollary~\ref{proof_of_order_relations}.

\begin{lemma}
\label{basic_results_about_Rich_class} Let the notation be as
described above.
\begin{enumerate}
\compress\romanate 

\item Contrary to our usual assumptions, let $G\le \GL(V)$ be any
algebraic group, let $C_\lambda$ and $C_\mu$ be any two unipotent
classes of $G$ with Jordan blocks given by the partitions $\lambda$
and $\mu$ respectively.  If $C_\lambda \le C_\mu$ then $\lambda \le
\mu$.

\item
Let $g\in Q$ and let $V_r$ be a subspace formed by $r$ Jordan blocks
of $g$.  Then $\dim V_r \le \sum_{i=1}^\ell \min\{r, \dim
W_i/W_{i-1}\}$ with equality holding if and only if $ \dim V_r\cap W_j
= \sum_{i=1}^j \min\{r, \dim W_i/W_{i-1}\} $ for all $j\ge 1$.
In particular, for all $r\ge1$ we have $\sum_{i=1}^r \lambda_i \le
\sum_{i=1}^\ell \min\{r, \dim W_i/W_{i-1}\}$.


\item Let $G\in \{\SO_{2n+1}, \SO_{2n}, \Sp_{2n}\}$.  If
$G=\SO_{2n+1}$ let $r=1$ and otherwise let $r=2$.  If there exists
$g\in Q$ with $(\lambda_1(g),\dots, \lambda_r(g)) = (\lambda_1, \dots,
\lambda_r)$ such that $\lambda_1(g), \dots, \lambda_r(g)$ are
nonsingular as Jordan blocks of $g$ then $\lambda_1, \dots, \lambda_r$
are nonsingular as Jordan blocks of $u$.
\end{enumerate}
\end{lemma}

\begin{proof}[Sketch of proof]
Part (i).  Let $U$ be the variety of all unipotent elements in $G$. 
It is easy to show that for each $j,b\ge 0$ the subset $\{g\in U \st
\dim \ker (g-1)^j \ge b\}$ is closed in $U$.  (One way to prove this
is to use elementary characterizations of rank in terms of
determinants of minors of a matrix.  Another way is to use the upper
semi-continuity of dimension applied to the endomorphism of $U\times
V$ given by $(g,v)\mapsto (g, (g-1)^jv)$, see
\cite[III.8.1]{spaltenstein}.)\ \ Let $u_\lambda^{}\in C_\lambda$ and
$u_\mu\in C_\mu$.  Since $C_\lambda \subseteq \overline {C_\mu}$ we
have that $u_\lambda$ is contained in any $G$-invariant, closed subset
of $U$ that contains $u_\mu$.  Thus, for each $j$, one has $u_\lambda
\in \{g\in U\st \dim \ker(g-1)^j \ge \dim \ker(u_\mu-1)^j\}$. 
Finally, note that $\dim \ker (u_\lambda-1)^j = \sum_{i=1}^j
\lambda_i^*$ and $\dim \ker (u_\mu -1)^j = \sum_{i=1}^j \mu_i^*$.

Part (ii) is elementary linear algebra and induction together with the
fact that $g$ acts trivially upon each factor $W_i/W_{i-1}$.


Part (iii) uses the following facts.  Every $P$-invariant, nonempty,
open subset of $Q$ contains $u$.  Let $X$ be any subset of $V$ and for
$g\in Q$ define a subspace $V_g:=\genby{(g-1)^i v \st i\ge 1, v\in X}
\le V$.  Then the set of $g\in Q$ such that $V_g$ is a nonsingular
subspace is an open set.  (Note that one can express the fact that
$V_g$ is nonsingular via a determinant being nonzero.)
\end{proof}

\begin{lemma}
\label{main_inductive_step}
If $G\in\{ \GL_n, \SO_{2n+1}\}$ let $r=1$, otherwise let $r=2$.  If
the following hypotheses hold then $\lambda = \mu$.
\begin{enumerate}
\compress\romanate

\item $(\lambda_1,\dots,\lambda_r) \le (\mu_1, \dots,\mu_r)$,

\item $\sum_{i=1}^\ell \min\{r,\dim W_i/W_{i-1}\} = 
\mu_1+\dots+\mu_r $,

\item there exists $g\in Q$ with $(\lambda_1(g), \dots, \lambda_r(g))
= (\mu_1,\dots, \mu_r)$ and if $G$ is orthogonal or symplectic 
$\lambda_1(g),\dots,\lambda_r(g)$ are nonsingular as Jordan block 
sizes of $g$.
\end{enumerate}
\end{lemma}

\begin{remark}
\label{remark_gather_pieces_for_main_proof}
The previous lemma abstracts the inductive step in showing $\lambda =
\mu$.  Essentially one can view Lemma~\ref{lemma_for_regular_elements}
and Lemma~\ref{lemma_giving_base_cases_for_two_blocks} as the base
cases.  Lemma~\ref{proof_of_main_thm_GLn} finishes the proof for the
case $G=\GL_n$.  For the remaining groups
Lemma~\ref{lemma_finding_mu_as_upper_bound} establishes part (i) and
(ii) and Lemma~\ref{lemma_giving_base_cases_for_two_blocks}
establishes part (iii).
\end{remark}

\begin{proof}
  If $G$ has rank $1$ then $P$ is a Borel subgroup and we are done by
  Lemma~\ref{lemma_for_regular_elements}.  We assume now that
  Theorem~\ref{theorem_Jordan_blocks_of_rich_classes} is true for
  classical groups with natural module $V'$ where $\dim V ' < \dim V$.
  
  Combining hypothesis (i) and
  Lemma~\ref{basic_results_about_Rich_class}(i) gives that $(\lambda_1
  (g)$, $\dots$, $\lambda_r(g))$ $\le$ $(\lambda_1$, $\dots$,
  $\lambda_r) \le (\mu_1$, $\dots$, $\mu_r)$, whence we have equality
  by hypothesis (iii).  Let $V_r$ be the space generated by $r$ Jordan
  chains of $u$, of lengths $(\lambda_1$, \dots, $\lambda_r) =
  (\mu_1$, \dots, $\mu_r)$.  If $G$ is symplectic or orthogonal we
  apply Lemma~\ref{basic_results_about_Rich_class}(iii) and assume
  that $V_r$ is nonsingular.  By hypothesis (ii) we have $\dim V_r =
  \lambda_1+\dots +\lambda_r = \mu_1+\dots + \mu_r = \sum_{i=1}^\ell
  \min\{r, \dim W_i/W_{i-1}\}$.

The inductive step will proceed as follows.  Let $X=V_r$.  We will
produce a $u$-stable decomposition $V=X\oplus Y$.  We will show that
$f$ induces flags in $X$ and $Y$ which we will denote by $f\cap X$ and
$f\cap Y$ such that $W_i$ is the direct sum of corresponding terms in
$f\cap X$ and $f\cap Y$.  We will then calculate the Jordan blocks of
the Richardson classes in $\Cl(X)$ and $\Cl(Y)$ (these are the
classical groups defined on $X$ and $Y$) associated with $f\cap X$ and
$f\cap Y$ and show that these equal $\lambda(u\restrict X)$ and
$\lambda(u\restrict Y)$.  The Jordan blocks of $u\restrict X$ are
$(\lambda_1,\dots,\lambda_r) = (\mu_1,\dots,\mu_r)$ by construction
(since $X=V_r$), and the blocks of $u\restrict Y$ will be found by
induction.

Let $f\cap X$ denote the flag in $X=V_r$ with terms given by
$X_i:=X\cap W_i$ for $1\le i \le \ell$.  We will construct below a
space $Y$ and a flag $f\cap Y$ with terms $Y_i$ such that $W_i = X_i
\oplus Y_i$ for $1\le i \le \ell$.  When $G$ equals $\Sp_n$ or $\SO_n$
the flags will be totally singular and in all cases $u$ will act
trivially upon the factors in each flag.

Let $G=\GL_n$.  Let $Y_1$ be any direct complement of $X_1$ in $W_1$;
this is $u$-stable since $u$ acts as $1$ on $W_1$.  Let $i\ge 2$ and
suppose $Y_{i-1}$ has been constructed such that $W_{i-1}=
X_{i-1}\oplus Y_{i-1}$.
%
Using the fact that $X_{i-1}\cap Y_{i-1}=\{0\}$ it is easy to show that
$\ker (u-1)\restrict{X_i}\cap \ker (u-1)\restrict{Y_{i-1}}=\{0\}$.
%
%
Then we may choose a direct complement $Z$ and a basis
$v_1,v_2,\dots$, of another direct complement as indicated:
$$\begin{array}{rcl} 
\ker(u-1)\restrict{W_i} &=& \ker(u-1)\restrict{X_i} \oplus 
\ker(u-1)\restrict {Y_{i-1}} \oplus Z\\
Y_{i-1}\cap (u-1) W_i & = & (u-1) Y_{i-1}\oplus \genby{v_1,v_2,\dots}.
\end{array}$$
%
%
%
For each $v_j$ fix $\hat v_j\in W_i$, a pre-image under $u-1$.  Let
$Y_i = Y_{i-1}\oplus \genby{\hat v_1, \hat v_2,\dots} \oplus Z$.  To
check that this sum is direct, write $0$ as the sum of an element in
each term, then apply $u-1$ and use the definitions.
%
%
It is now relatively easy to show that $W_i = X_i \oplus Y_i$.  
Finally, we take $Y=Y_\ell$.  

If $G\ne \GL_n$ let $Y= V_r^\perp$ and $Y_i = Y\cap W_i$ for $1\le i
\le \ell$.  The dimension of $X_i$ can be calculated using
Lemma~\ref{basic_results_about_Rich_class}(ii) and $\dim Y_i$ is given by
$\dim W_i +\dim(W_i^\perp \cap X) -\dim X$.  Using dimension
calculations one may show that $W_i=X_i \oplus Y_i$ for $1\le i \le
\ell$ and that the flags are totally singular.

For $J\in \{X,Y\}$ let $\Cl(J)$ be the classical group on $J$, let
$P_J$ be the parabolic in $\Cl(J)$ corresponding to the flag $f\cap J$
and let $Q_J$ be the unipotent radical of $P_J$.  We may identify
$Q_XQ_Y$ as a subgroup of $Q$.  Let $C$ denote the Richardson orbit in
$Q$ and note that $C\cap (Q_XQ_Y)$ is an open subset of $Q_XQ_Y$ which
is also dense as it contains $u$.  Then $C$ contains the Richardson 
orbits in $Q_X$ and $Q_Y$.  Let $u'\in C\cap (Q_XQ_Y)$ such that 
$u'\restrict X$ and $u'\restrict Y$ represent the Richardson orbits 
in $Q_X$ and $Q_Y$ respectively.  

We have that $u$ and $u'$ are conjugate whence $\lambda(u) =
\lambda(u')$.  We also have that $\lambda( u \restrict X) =
(\lambda_1,\dots, \lambda_r) = (\mu_1, \dots, \mu_r)$ by construction. 
We have $\lambda(u'\restrict X) \ge \lambda(u\restrict X)$ since
$u'\restrict X$ represents the Richardson orbit (and using
Lemma~\ref{basic_results_about_Rich_class}(i)\,).  This, together with the 
fact that $\lambda(u') = \lambda(u)$ implies that $\lambda(u'\restrict 
X) = \lambda(u\restrict X)$.  Thus $\lambda(u'\restrict Y) = 
\lambda(u\restrict Y)$ and we may assume, for our purposes, that 
$u=u'$.  

We have $\lambda(u\restrict X)= (\lambda_1,\dots,\lambda_r)$ and
$\lambda( u\restrict Y) = (\lambda_{r+1},\dots)$.  Since $(\lambda_1$,
$\dots$, $\lambda_r)$ $=$ $(\mu_1$, $\dots$, $\mu_r)$ it suffices to
show that $(\lambda_{r+1} ,\dots ) = (\mu_{r+1}, \dots )$.  One may
calculate the parts of the Levi partition of $Y$ using the dimensions
of factors in the flag $f\cap Y$.
%
%
One may verify that $\Lambda(Y)=(\max\{n_i-r, 0\} \st 1\le i \le s)$
if $m=0$ and $\Lambda(Y)= (\max\{n_i-r,0\} \st 1\le i \le s) \oplus
(m-1)$ if $m\ge 1$ and $G\ne \SO_{2n+1}$ and $\Lambda(Y) =
(\max\{n_i-1,0\} \st 1\le i \le s)\oplus m$ if $G=\SO_{2n+1}$ (when
$G=\SO_{2n+1}$ then $m$ in $\Lambda$ corresponds to $\SO_{2m+1}$, but
$Y$ is even dimensional and $m-1$ or $m$ in $\Lambda(Y)$ corresponds
to $\SO_{2m-2}$ or $\SO_{2m}$).
%
%
%
%
%

By induction we may apply
Theorem~\ref{theorem_Jordan_blocks_of_rich_classes} to determine the
Jordan blocks of this Levi partition $\Lambda(Y)$.  By analyzing the
cases in Theorem~\ref{theorem_Jordan_blocks_of_rich_classes}, one
finds that they equal $\mu$ with the first $r$ rows removed.
%
%
%
%
%
%
\end{proof}

\begin{numproof}[Proof of
Theorem~\ref{theorem_Jordan_blocks_of_rich_classes} when $G=\GL_n$]  
\label{proof_of_main_thm_GLn}
Note that $\mu_1 = \ell = s$.  Using
Lemma~\ref{basic_results_about_Rich_class}(ii) it is easy to verify
hypotheses (i) and (ii) of Lemma~\ref{main_inductive_step}.

It remains to prove the existence of $g\in Q$ with $\lambda_1(g)=
\mu_1$.  Let $X_1$ be a one dimensional subspace of $W_1$ and $Y_1\le
W_1$ such that $W_1=X_1\oplus Y_1$.  For $i\in \{2,\dots,\ell\}$ let
$X_i$ be an $i$ dimensional subspace of $W_i$ such that $X_{i-1} \le
X_i$ and let $Y_i\le W_i$ such that $Y_{i-1}\le Y_i$ and $W_i = X_i
\oplus Y_i$.  Define $f\cap X$ to be the flag in $X=X_\ell$ with terms
given by the $X_i$ and $f\cap Y$ to be the flag in $Y=Y_\ell$ with
terms given by the $Y_i$.  Let $P_X$ be the parabolic in $\GL(X)$ of
$f\cap X$, let $Q_X$ be the unipotent radical of $P_X$ and identify
$Q_X$ as a subgroup of $Q$.  Then $P_X$ is a Borel subgroup of
$\GL(X)$, whence there exists an element $g$ in $Q_X$ which has one
block of size $\ell=\mu_1$ by Lemma~\ref{lemma_for_regular_elements}.
\end{numproof}

\begin{cor} \label{every_class_in_GLn_is_richardson}%
Every unipotent class in $\GL_n$ is a Richardson class.
\end{cor}

This is also proven in \cite[II.5.14]{spaltenstein} and in \cite[5.5]{humphreys3}.  

\begin{proof}
If a unipotent class has Jordan blocks given by the partition $\nu$
then it is the Richardson class of any parabolic with Levi partition
equal to $\nu^*$.
\end{proof}

For the next result we introduce some notation.  Let $H\le J$ be
algebraic groups and let $H$ act upon $J$ via conjugation.  Given a
subset $O\subseteq J$ we denote by $\overline O$ the closure taken
within $J$ and by $O^J$ the subset $\bigcup_{g\in J} O^g = \{gxg\inv
\st g\in J, x\in O\}$.

\begin{lemma}
\label{lemma_for_descent_in_partial_order}
Let $H\le J$ be algebraic groups and use the notation described
above.  Let $O_1$ and $O_2$ be two $H$-classes in $H$.  
\begin{enumerate}
\compress\romanate 

\item If $O_1 \subseteq \overline{O_2}$ then $O_1^J \subseteq
\overline{O_2^J}$.

\item If $H$ has a dense orbit in $H\cap \overline{O_2^J}$, has a
single orbit in $H\cap O_2^J$ (i.e. $H\cap O_2^J = O_2$), and $O_1^J$
is a subset of $\overline{O_2^J}$ then $O_1 \subseteq \overline{O_2}$.

\item If $H$ has a single orbit in $O_1^J$ (i.e. $O_1^J = O_1$), has
finitely many orbits in $O_2^J$, and $O_1^J$ is a subset of $\overline
{O_2^J}$ then $O_1 \subseteq \overline{O_2}$.

\end{enumerate}
\end{lemma}

We refer to conditions (ii) and (iii) as ``descending from $J$ to
$H$''.

\begin{proof}
Part (i) We have:
$$O_1^J \subseteq \big(\,\overline{O_2}\,\big)^J = \bigcup_{g\in J}
\big(\,\overline{O_2}\,\big)^g = \bigcup_{g\in J} \overline{O_2^g} \subseteq
\overline{\bigcup_{g\in J} O_2^g} = \overline{O_2^J}.$$

Part (ii) We claim that: 
$$O_1 \subseteq H\cap O_1^J \subseteq H\cap \overline{O_2^J} = 
\overline{O_2}.$$
The final equality is the one to be proved.  Let $C$ denote a dense
orbit of $H$ in $H\cap \overline{O_2^J}$.  Then $C\subseteq \overline
{O_2^J}$ whence $C^J \subseteq \overline{O_2^J}$ by part (i).  On the
other hand, $O_2\subseteq \overline C$ whence $O_2^J \subseteq
\overline{C^J}$ by (i).  Thus $O_2^J = C^J$ whence $C\subseteq
H\cap C^J = H \cap O_2^J = O_2$ and $C=O_2$.

Part (iii).  Denote the $H$-orbits in $O_2^J$ by $O_2=O_2^{g_1}$,
$O_2^{g_2}$, \dots, $O_2^{g_b}$ where $g_1=1$ and $g_i\in J-H$ for
$i>1$.  We have: $O_1 \subseteq \overline{O_2^J} =
\overline{O_2^{g_1}} \cup \dots \cup \overline{O_2^{g_b}}$ which
implies that $O_1\subseteq \overline{O_2^{g_j}}$ for some $j$.  Then
$O_1 = O_1^{g_j\inv} \subseteq \Big(\overline{O_2^{g_j}}\Big)^{g_j\inv} = 
\overline {O_2}$.

\end{proof}

\begin{cor}\label{proof_of_order_relations}
Contrary to our usual assumptions, let $G\le \GL_n$ be an algebraic
group, let $C_\lambda$ and $C_\mu$ be two unipotent classes with 
Jordan blocks given by the partitions $\lambda$ and $\mu$.

\begin{enumerate}\compress\romanate

\item Suppose that all the unipotent elements in $G$ with Jordan
blocks equal to $\mu$ form a single conjugacy class.  Then $\lambda <
\mu$ if and only if $C_\lambda < C_\mu$.

\item If $G=\GL_n$, or $G\in \{\O_n,\Sp_n\}$ and $p\ne2$, or $G\in
\{\O_n, \Sp_n\}$ and $\mu$ has no even parts with even multiplicity,
or $G= \SO_{n}$ with $n$ even and $p\ne 2$, then $\lambda < \mu$ if and
only if $C_\lambda < C_\mu$.

\end{enumerate}
\end{cor}

\begin{proof}
Part (i).  By Lemma~\ref{basic_results_about_Rich_class}(i) we have
that $C_\lambda < C_\mu$ implies $\lambda < \mu$.  We prove the
converse first for $G=\GL_n$.

Step 1: Since $\lambda \le \mu$ we may fix a sequence of partitions
$\lambda = \lambda^{(0)} < \lambda^{(1)} < \dots < \lambda^{(r)} =
\mu$ such that for each $i$ we have that $\lambda^{(i)}$ and
$\lambda^{(i+1)}$ differ in exactly two places, i.e. there exist
exactly two indices $j$ such that $\lambda_j^{(i)} \ne
\lambda_j^{(i+1)}$ (see \cite[p23]{james_kerber}).  Then by
transitivity it suffices to prove that $\lambda < \mu \implies
C_\lambda < C_\mu$ when $\lambda$ and $\mu$ differ in exactly two
places, which we now assume.

Step 2: Since $\lambda$ and $\mu$ differ in exactly two places, we may
find a subgroup $\GL_{n_1}\GL_{n_2}$ of $G$, a unipotent
$\GL_{n_1}$-class $C$, two unipotent $\GL_{n_2}$-classes $C_1$ and
$C_2$ with $C_\lambda$ and $C_\mu$ the extensions to $G$ of $CC_1$ and
$CC_2$ respectively (i.e. the classes $C_1$ and $C_2$ correspond to
the two parts where $\lambda$ and $\mu$ differ).  By
Lemma~\ref{lemma_for_descent_in_partial_order}(i) it suffices to show
that $C_1\le C_2$ (for then $CC_1 \le CC_2$ and $C_\lambda = (CC_1)^G
\le (CC_2)^G = C_\mu$).

Step 3.  It suffices now to prove the result under the assumption that
$\lambda$ is a two part partition (whence $\mu$ has one or two parts).
Then the difference between $\mu^*$ and $\lambda^*$ is that
$\lambda^*$ has one extra $2$ and two fewer $1$'s.  Let $g\in
C_\lambda$.  By Corollary~\ref{every_class_in_GLn_is_richardson}, we
may find flags $f_\lambda: 0< W_2 < W_3< \dots$ and $f_\mu: 0 < W_1 <
W_2 < W_3 < \dots$ such that $f_\lambda$ and $f_\mu$ have
corresponding Levi partitions of $\lambda^*$ and $\mu^*$, $f_\lambda$
and $f_\mu$ are identical to the right of $W_3$, and $g$ represents
the Richardson orbit corresponding to $f_\lambda$ (in particular $g$
acts trivially upon each factor in $f_\lambda$).  Then $g$ is in the
unipotent radical associated with $f_\mu$, which in turn is contained
in $\overline {C_\mu}$.  Whence, $C_\lambda \subseteq
\overline{C_\mu}$.

Now (i) is proven for $G=\GL_n$.  If $G<\GL_n$ we may descend to $G$
via part (ii) of the previous lemma; i.e., apply
Lemma~\ref{lemma_for_descent_in_partial_order}(ii) with $H=G$ and
$J=\GL_n$ to get $\lambda<\mu$ implies $C_\lambda < C_\mu$.

Part (ii).  This is immediate from part (i) (and the comments in
Section~\ref{conventions}), unless $G=\SO_{n}$ with $n$ even, $p\ne2$. 
However, part (i) holds for $\O_n$ and one may descend to $\SO_n$ by
applying Lemma~\ref{lemma_for_descent_in_partial_order}(iii) with
$H=\SO_n$ and $J=\O_n$.
\end{proof}

The following lemma establishes
Lemma~\ref{main_inductive_step}(i),(ii) for the cases where $G\ne
\GL_n$.  It will also be used in
Lemma~\ref{lemma_giving_base_cases_for_two_blocks} to establish
Lemma~\ref{main_inductive_step}(iii).

\begin{lemma}
\label{lemma_finding_mu_as_upper_bound}
If $G=\SO_{2n+1}$ let $r=1$ and if $G\in \{\Sp_{2n}, \SO_{2n}\}$ let
$r=2$.  Recall that $\ell$ is the number of terms in the natural flag
and that the Levi partition is $\Lambda = (1^{c(1)}, 2^{c(2)}, \dots)
\oplus m$.  Then $(\mu_1,\dots,\mu_r)$ are listed below.  Furthermore,
$(\lambda_1, \dots,\lambda_r) \le (\mu_1,\dots,\mu_r)$ and
$\sum_{i=1}^\ell \min\{r,\dim W_i/W_{i-1}\} = \mu_1+\dots+\mu_r$.
\begin{enumerate}
\romanate\compress

\item $G=\SO_{2n+1}$ and $p\ne 2$.  We have $\mu_1 = \ell$.  

\item $G=\SO_{2n}$ and $p\ne 2$.  If $m=0$ and $c(1)\ge 1$ then 
$(\mu_1,\mu_2) = (\ell - 1, \ell - 2c(1) +1)$.  Otherwise we 
have $(\mu_1,\mu_2) = (\ell, \ell -2c(1))$.

\item $G=\SO_{2n}$ and $p=2$.  If $m=0$ and $c(1)\ge 2$ then
$(\mu_1,\mu_2) = (\ell -2, \ell -2c(1)+2)$.  If $m=0$, $c(1)=1$, or
$m\ge 1$, $c(1)\ge 1$ then $(\mu_1,\mu_2) = (\ell -1, \ell-2c(1)+1)$.  
If $c(1)=0$ then $(\mu_1,\mu_2) = (\ell, \ell -2c(1))$.  

\item $G=\Sp_{2n}$.  If $m\ge1$ and $c(1) \ge 1$ then
$(\mu_1,\mu_2) = (\ell - 1, \ell -2c(1)+1)$.  Otherwise
$(\mu_1,\mu_2) = (\ell, \ell -2c(1))$.
\end{enumerate}
\end{lemma}

%

\begin{proof}
  Recall that $\psi (\Lambda)$ is defined in
  table~\ref{table_of_Jordan_blocks_of_Rich_class} and that $\mu$
  equals $\psi(\Lambda)^*$, the dual of $\psi(\Lambda$).  Thus $\mu_1$
  equals the number of parts in $\psi(\Lambda)$ and $\mu_2$ equals the
  number of parts in $\psi(\Lambda)$ which are greater than or equal
  to $2$.  It is easy to verify the stated formulas for $\mu_1$ and
  $\mu_2$.

Since $\sum_{i=1}^\ell \min\{1, \dim W_i/W_{i-1}\}=\ell$ and
$\sum_{i=1}^\ell \min\{2,\dim W_i/W_{i-1}\} = 2\ell -2c(1)$ we
conclude that $\sum_{i=1}^\ell \min\{r, \dim W_i/W_{i-1}\} = \mu_1 +
\dots + \mu_r$ and that $(\lambda_1,\lambda_2) \le (\ell, \ell-2c(1))$
%
%
(by Lemma~\ref{basic_results_about_Rich_class}(ii)\,).  
This gives the desired upper bounds on $\lambda$ in case (i) or when
$(\mu_1,\mu_2) = (\ell, \ell -2c(1))$.
%

For all the remaining cases we start by proving that one cannot have
$\lambda_1 = \ell$, whence $(\lambda_1,\lambda_2) \le (\ell -1,
\ell-2c(1)+1)$.  First we assume $G=\SO_{2n}$, $m=0$ (whence $\ell =
2s$) and $c(1)\ge 1$.  If $\lambda_1 = \ell$ there is a Jordan chain
$v_1$,\dots,$v_{2s}$ of length $2s=\ell$.  Since $u$ acts trivially
upon each factor in the flag, and since the Jordan chain has as many
elements as there are terms in the flag, we see that $v_s\in
W_{s}-W_{s-1}$ and $v_{s+1}\in W_{s+1} - W_s$.  By
Lemma~\ref{jordan_chain_degeneracy}(i) we have $\beta(v_s,v_{s+1})\ne
0$.  Let $\tilde W= \genby{v_{s+1}}^\perp \cap W_{s}$.  Then $\tilde
W$ is a totally singular $(n-1)$-space, whence $\tilde W^\perp/\tilde
W$ is a nonsingular $2$-space.  But $v_s$ and $v_{s+1}$ project to
distinct, nontrivial elements in $\tilde W^\perp /\tilde W$ whence
$u\restrict {\tilde W^\perp/\tilde W}$ is a nontrivial unipotent
element of $\SO(\tilde W^\perp/\tilde W) = \SO(2)$, a contradiction.
In all other cases where $(\mu_1,\mu_2) < (\ell, \ell-2c(1))$ the
parity conditions upon $\lambda$ and the fact that $\lambda_1 +
\lambda_2 \le 2\ell -2c(1)$ imply $\lambda_1\ne \ell$.

It remains to show that if $G=\SO_{2n}$, $m=0$ and $c(1)\ge 2$ then
$(\lambda_1,\lambda_2) \le (\ell-2, \ell-2c(1)+2)$.  However, we have
shown already that $(\lambda_1,\lambda_2)\le (\ell-1,\ell-2c(1)+1)$
and if $\lambda_1=\ell-1$ this contradicts the parity conditions upon
$\lambda$.
\end{proof}

\begin{lemma} 
\label{lemma_giving_base_cases_for_two_blocks}
With the usual notation, the following hold and, in particular,
Lemma~{\rm\ref{main_inductive_step}}(iii) holds, whence 
Theorem~{\rm\ref{theorem_Jordan_blocks_of_rich_classes}} is proven.

\begin{enumerate}
\compress\romanate

\item Let $G=\SO_{2n}$ and $p=2$.  {\rm (i)(a)} If $\Lambda = (2^a)$
then $\lambda = (2a,2a)$.  {\rm (i)(b)} If $\Lambda=(2^a,1)$ then
$\lambda=(2a+1,2a+1)$.  {\rm (i)(c)} If $\Lambda = (2^a,1^b)$ with
$b\ge 2$ then $\lambda=(2a+2b-2,2a+2)$.

\item Let $G= \SO_{2n}$ and $p\ne 2$.  {\rm (ii)(a)} If $\Lambda =
(2^a)$ then $\lambda= (2a,2a)$.  {\rm (ii)(b)} If $\Lambda= (2^a,1^b)$ 
with $b\ge 1$ then $\lambda = (2a+2b-1,2a+1)$.

\item Let $G=\Sp_{2n}$.  {\rm (iii)(a)} If $\Lambda = (2^a)$ then
$\lambda = (2a,2a)$.  {\rm (iii)(b)} If $\Lambda = (2^a,1^b)$ then
$\lambda = (2a+2b,2a)$.  {\rm (iii)(c)} If $\Lambda = (2^a)\oplus 1$
then $\lambda = (2a+1,2a+1)$.  {\rm (iii)(d)} If $\Lambda = (2^a,1^b)
\oplus 1$ with $b\ge 1$ then $\lambda = (2a+2b, 2a+2)$.
\end{enumerate}
\end{lemma}

We sketch two proofs of parts (i)-(iii).  Neither proof seems entirely
satisfactory as each contains a tedious verification of a rather
simple fact.

\begin{proof}
\emph{Sketch of first proof of {\rm (i)-(iii)}}.  For the following
statement, $\lambda$ and $\mu$ need not have their usual definitions. 
Let $\lambda$ and $\mu$ be two partitions and $u_\lambda$, $u_\mu$ two
unipotent elements with Jordan blocks given by $\lambda$ and $\mu$
respectively.  Suppose $\mu$ has two parts.  Then we claim that
$\lambda < \mu$ implies $\dim \C G(u_\lambda) \ge \dim \C G(u_\mu)$
with the inequality strict provided the parts of $\mu$ are
nonsingular.

Given the claim, we again let $\lambda$ and $\mu$ have their usual
definitions, whence $\lambda \le \mu$ by
Lemma~\ref{lemma_finding_mu_as_upper_bound}.  Let $L$ be the Levi
factor of the parabolic $P$ under discussion.  Then Richardson's
Theorem (iii) and the claim show $\dim L = \dim \C G(u_\lambda) \ge
\dim \C G(u_\mu)$.  Using the expression for $\mu$ in the previous
Lemma, one may check that $\dim \C G(u_\mu)$ equals $\dim L$ or $\dim
L+2$ with the latter only if both parts of $\mu$ are singular.  We
conclude that both parts of $\mu$ are nonsingular, and that we must
have $\lambda = \mu$.

(We note that it would be circular to prove the claim by applying
Spaltenstein's version of Theorem~\ref{theorem_on_order_relations}, as
this is proven in \cite{spaltenstein} using the version of
Theorem~\ref{theorem_Jordan_blocks_of_rich_classes} which is found
there.)  To prove the claim directly one may manipulate the formulas
for dimensions of centralizers, though this is somewhat tedious.  In
particular, let $c$ be the multiplicity of $\lambda_1^*$ in
$\lambda^*$.  If $\lambda_1^*=5$, or $\lambda_1^*=4$, $c\ge2$, or
$\lambda_i^* =3$, $c \ge 4$ then one may show that $\sum_{i\ge
  1}\big((\lambda_i^*)^2-(\mu_i^*)^2\big) > 2\lambda_1^*$ which proves
the claim (by examining the formulas for dimensions of centralizers).
The remaining cases amount to direct calculations.

\emph{Sketch of second proof of {\rm (i)-(iii)}.}  There are two cases.

Case 1: $\mu$ has two equal parts.  We claim that there exists $g\in
Q$ with $\mu =\lambda(g)$.  Given the claim, and using
Lemma~\ref{basic_results_about_Rich_class}(i) and
Lemma~\ref{lemma_finding_mu_as_upper_bound}, we have $\mu = \lambda(g)
\le \lambda \le \mu$.

To prove the claim let $\mu = (n,n)$ where $n$ is the rank of $G$. 
The tedious part of the argument is verifying, inductively, that one
may construct, using roots in $\rtsys(Q)$, a root base of an $A_{n-1}$
root system.  Given this root base, the group generated by the maximal
torus and the root groups corresponding to $\Z$-linear combinations of
this base is isomorphic to $\GL_n$.  Let $g$ in $\GL_n$ be a regular
unipotent element written as the product of a nontrivial element in
each root group corresponding to a root in this root base (see
\cite{steinberg_regular_elements}).  Then $g$ is in $Q$ and $g$ has
two blocks of size $n$ in the natural embedding of $\GL_n$ in $G$.

Case 2: $\mu$ has two distinct parts.  As stated in
Section~\ref{conventions}, one sees that $G$ has a unique unipotent
class $C_\mu$ with Jordan blocks given by $\mu$.  Let $C_\lambda$ be
the Richardson class of $P$.  By
Lemma~\ref{lemma_finding_mu_as_upper_bound} we have $\lambda \le \mu$,
and by Lemma~\ref{proof_of_order_relations} we have $C_\lambda \le
C_\mu$.  One may easily show that $\dim C_\mu = \dim C_\lambda$ whence
$C_\lambda = C_\mu$ and $\lambda = \mu$.

\emph{Sketch of proof of {\rm \ref{main_inductive_step}(iii)}.} This
proof parallels that given for the case $G=\GL_n$ (see
Proof~\ref{proof_of_main_thm_GLn}) so we will be brief.  Recall that
the natural flag $f$ has terms $W_i$.  We will produce a decomposition
of $f$, by constructing and decomposing a flag $\tilde f$ which is
isomorphic to (whence conjugate to) $f$.

Suppose $G= \SO_{2n+1}$.  Let $\tilde X\le V$ be a nonsingular
subspace of dimension $2s+1=\ell$.  Choose a totally singular flag
$0<\tilde X_1 < \dots < \tilde X_\ell = X$ where $\dim \tilde X_i =
i$.  This flag corresponds to a Borel subgroup of $\Cl(\tilde X)$
which has Levi partition $\Lambda(\tilde X) = (1^s)$.  Suppose
$G\in\{\Sp_{2n}, \SO_{2n}\}$.  Let $\tilde X\le V$ be a nonsingular
subspace of dimension $2\ell-2c(1)$.  Choose a totally singular flag
$0<\tilde X_1 < \dots < \tilde X_\ell = X$ where $\dim X_j =
\sum_{i=1}^j \min\{2, \dim W_i/W_{i-1}\}$.  This flag has Levi
partition $\Lambda(\tilde X)$ as follows: If $m=0$ then
$\Lambda(\tilde X) = (2^{s-c(1)}, 1^{c(1)})$; If $m\ge 1$ and
$G=\Sp_{2n}$ then $\Lambda(\tilde X) = (2^{s-c(1)}, 1^{c(1)})\oplus
1$; If $m\ge 1$ and $G=\SO_{2n}$ then $(2^{s-c(1)}, 1^{c(1)+1})$.

In each case we define $\tilde Y= \tilde X^\perp$ and choose a totally
singular flag $0 \le \tilde Y_1 \le \dots \le \tilde Y_\ell = \tilde
Y$ where $\dim Y_i = \dim W_i - \dim X_i$.  Define the flag $\tilde f$
to have terms $0< \tilde W_1 < \dots < \tilde W_\ell = V$ where
$\tilde W_i = \tilde X_i \oplus \tilde Y_i$.  Since $\tilde f$ is
conjugate to $f$ we see that a similar decomposition holds for $f$
which we express as $V=X\oplus Y$, $f= (f\cap X)\oplus (f\cap Y)$. 
The Levi partitions for the flag $f\cap X$ are the Levi partitions
listed for $\Lambda(\tilde X)$ above.

Now that $f$ has been decomposed, let $P_X$ be the parabolic in
$\Cl(X)$ and let $Q_X$ be its unipotent radical.  Then we identify
$Q_X$ as a subgroup of $Q$.  Let $g\in Q_X\le Q$ which represents the
Richardson orbit in $Q_X$.  Then we may apply parts (i)-(iii) to
calculate the Jordan blocks of $g$.  We find that $\lambda(g) = \mu_1$
if $G=\SO_{2n+1}$ and $\lambda(g) = (\mu_1,\mu_2)$ if $G\in
\{\SO_{2n}, \Sp_{2n}\}$ where $\mu_1$ and $\mu_2$ are as in
Lemma~\ref{lemma_finding_mu_as_upper_bound}.
%
%
%
\end{proof}

\section{Richardson Classes of Distinguished Parabolics}

\begin{lemma}
\label{partitions_of_distinguished_richardson_classes} Let $G$ be one
of $\GL_n$, $\SO_{2n+1}$, $\SO_{2n}$, and $\Sp_{2n}$.  Let $\Psi$
denote the map taking each distinguished parabolic class to the Jordan
blocks of its Richardson class.  Then $\Psi$ gives a bijection with
the set of partitions described in
table~{\rm\ref{table_jordan_blocks_dist_rich_classes}}. 
\end{lemma}

For $p\ne 2$, the descriptions in
table~\ref{table_jordan_blocks_dist_rich_classes} of the image of
$\Psi$ are stated in \cite{bala-carter}, but it is not stated there
that these partitions equal the Jordan blocks of the Richardson class.

\begin{table*}
\caption{Jordan blocks of distinguished Richardson classes}%
\label{table_jordan_blocks_dist_rich_classes}
$$\begin{array}{|p{4in}|}
\hline
\hspace{1in}Image of\ \ $\Psi$\\
\hline
\qquad $\GL_n$\\
The partition of $n$ consisting of a single block\\
\hline
\qquad $G= \SO_{2n+1},\ p\ne 2$\\
Partitions of $2n+1$ consisting of distinct 
odd parts\\
\hline
\qquad $G=\SO_{2n+1},\ p=2$\\
Partitions of $2n+1$ of the form $1\oplus
\lambda$ such that: each part of $\lambda$ is even; the 
multiplicity of each part of $\lambda$ is at most $2$; if $i$ 
is even then $\lambda_i - \lambda_{i+1} \ge 4$.\\
\hline
\qquad $G=\Sp_{2n}$\\
 Partitions of $2n$ consisting of distinct even 
parts\\
\hline
\qquad $G= \SO_{2n},\ p\ne 2$\\
 Partitions of $2n$ consisting of distinct odd 
parts\\
\hline
\qquad $G=\SO_{2n},\ p=2$\\
Partitions $\lambda$ of $2n$ such that: $\lambda$
has an even number of parts; each part of $\lambda$ is even; the
multiplicity of each part is at most $2$; if $i$ is even with
$\lambda_{i+1}\ne 0$ then $\lambda_{i}-\lambda_{i+1} \ge 4$.\\
\hline
\end{array}$$
\end{table*}

\begin{proof}
If $G=\GL_n$, then the only distinguished parabolic is the Borel 
subgroup, which corresponds to the regular class. 

We give the proof for $\SO_{2n}$ and leave the other cases (which are
simpler) to the reader.

Let $\Lambda$ be the Levi partition of a distinguished parabolic $P$. 
Using the description of distinguished parabolics given in
\cite{carter} we may write $\Lambda =(n_1, \dots, n_s) \oplus m =
(1^{c(1)}, \dots, (2m)^{c(2m)}, (2m{+}1)^{c(2m{+}1)})\oplus m$ where
we index the $n_i$ such that $n_{s}$ is the largest $n_i$.  If
$m=0$ then $n_{s}\in \{1,2\}$; if $m\ge 1$ then  $n_{s} \in
\{2m-1, 2m\}$ and in all cases $c(i) \ge 1$ if and only if $1\le i
\le n_s$.

Let $\psi$ be the map defined in
Theorem~\ref{theorem_Jordan_blocks_of_rich_classes}.  If $p\ne2$ and
$m=0$ then $\psi (\Lambda) = (1^{2c(1)-2},2^{2c(2)+1})$.  If $p\ne 2$
and $m\ge2$ then $\psi (\Lambda) = (1^{2c(1)}, \dots,
(2m-1)^{2c(2m-1)}, (2m)^{c(2m)+1})$.  If $p=2$ and $m=0$ then
$\psi(\Lambda) = (1^{2c(1)-4},2^{2c(2)+2})$.  If $p=2$ and $m\ge 2$
then $\psi(\Lambda) = (1^{2c(1)-2}$, $2^{2c(2)+2}$, $\dots$,
$(2m-1)^{2c(2m-1)-2}$, $(2m)^{2c(2m)+2})$.


We have that $\Psi(P)$ equals $\psi(\Lambda)^*$, the dual of
$\psi(\Lambda)$.  The formulas for $\psi(\Lambda)$ make it clear that
$\Psi$ is injective.

Let $\mu$ be any partition and let $m(i)$ be the multiplicity of $i$
in the dual partition $\mu^*$.  Then $\mu$ consists of distinct parts
if and only if $\mu^*$ contains each integer between $1$ and its
maximal part.  Also, each part of $\mu$ is odd if and only if for each
$j$ we have $\sum_{i\ge j} m(i)$ is odd.  

If $P$ is given and $p\ne 2$, the comments just made about $\mu$ show
that $\Psi(P) = \psi(\Lambda)^*$ satisfies the properties described in
the statement of the Lemma.  In other words, the image of $\Psi$ is in
the desired set.

Conversely, let $p\ne2$ and let $\lambda$ be given which satisfies the
properties described in the statement of the Lemma.  The comments just
made about $\mu$ show that $\lambda^*$ can be set equal to an
expression of the form given for $\psi(\Lambda)$, and then one may
solve for $c(1)$, $c(2)$, etc.  In other words, $\Psi$ is surjective.

The case for $p=2$ may be verified similarly, however the following
alternative description may make the proof easier.  Let $\lambda\two$
and $\lambda\nottwo$ be the Jordan blocks of the Richardson class of a
parabolic associated with $\Lambda$ when $p=2$ and when $p\ne 2$
respectively.  Let $\lambda\nottwo = (\lambda_1, \lambda_2, \dots,
\lambda_{2\ell-1}, \lambda_{2\ell})$ where $\lambda_{2\ell-1}$ or
$\lambda_{2\ell}$ is the last nonzero part of $\lambda\nottwo$.  Then
$\lambda = (\lambda_1-1, \lambda_2+1, \dots, \lambda_{2\ell-1}-1,
\lambda_{2\ell}+1)$.  This description may be verified directly from
the formulas for $\psi(\Lambda)$ (Spaltenstein
\cite[III.7.2,III.8.2]{spaltenstein} defines a similar map for the
Jordan blocks of all unipotent classes; note there is a typographical
mistake in the formula for $\SO_{2n+1}$).
%

Given $\lambda = \psi(\Lambda)^*$, it remains to prove that $\lambda$
is nonsingular.  For those cases where $\lambda$ has distinct parts
this follows from Lemma~\ref{jordan_chain_degeneracy}.  In the
remaining cases we have that $G$ is orthogonal and $p=2$.  By 
Richardson's Theorem (iii) we know $\dim L = \dim \C G(u)$ where
$L$ is the Levi subgroup determined by $\Lambda$ and $u$ is an element
of the Richardson class in $G$.  It is now easy to finish the proof 
by using Spaltenstein's expression for $\dim \C G(u)$ described 
in Remarks~\ref{remark_on_spaltenstein_centralizer_formula}.  
%
%
\end{proof}

\begin{cor}
Let \label{richardson_map_inj_on_dist} $G$ be a simple algebraic group
and consider the map which takes each distinguished parabolic class to
its Richardson class.  This map is injective.
\end{cor}

\begin{proof}
  For the classical groups this follows from the previous lemma.  For
  the exceptional groups, we observe that no two distinct
  distinguished parabolics have the same dimension of Levi factor.  By
  Theorem~\ref{Richardsons_theorem} (iii) the dimension of the Levi
  factor equals the dimension of the centralizer of an element in the
  unipotent class, whence the result follows by dimension.
\end{proof}

\section{Proof of the Bala-Carter-Pommerening Theorem}
Throughout this section, $G$ denotes a connected reductive group, 
unless indicated otherwise.  

\begin{lemma}
\label{min_Levi_unique_up_to_conj}
\begin{enumerate}
\romanate\compress \item Let $S$ be a torus in $G$.  Then $L=\C G(S)$
is a Levi subgroup.

\item If $u$ is a unipotent element and $S$ a maximal torus of $\C
G(u)$ then $u$ is distinguished in $L=\C G(S)$.  Furthermore, any Levi
subgroup in which $u$ is distinguished is conjugate to $L$ via an
element of $\C G(u)^\circ$.
\end{enumerate}
\end{lemma}

\begin{proof}  For part (i) one may adapt \cite[5.9.2]{carter}.  For 
part (ii) one may adapt \cite[5.9.3]{carter}.
\end{proof}

%
%

\begin{cor}
\label{bijection_reduction_to_dist_classes} Define a map from
$G$-classes of pairs $(L,C)$ consisting of a Levi subgroup $L$ of $G$
and a distinguished unipotent $L$-class $C$ to unipotent $G$-classes
by extending $C$.  This map gives a bijection.
\end{cor}

%

%

\begin{lemma}
\label{rep_of_dense_orbit_dist}
Let $P$ be a distinguished parabolic of $G$.  Let $\overline G=
G/Z(G)$, $\overline P = P/Z(G)$, let $\overline Q$ be the unipotent
radical of $\overline P$ and let $u$ represent the dense orbit of
$\overline P$ upon its unipotent radical $\overline Q$.  Then $\C
{\overline G}(u)^\circ = \C {\overline P}(u)^\circ = \C {\overline
Q}(u)^\circ$.  In particular the Richardson class of $P$ is
distinguished in $G$.
\end{lemma}

\begin{proof}
It is easy to reduce to the case $Z(G)=1$ and adapt the proof given 
in \cite[5.8.7]{carter}.  
\end{proof}


\begin{numproof}[Proof of Theorem~\ref{bala-carter-theorem}]
\label{proof_of_bala_carter}
Part (i). This is by definition of the map $\psi$.

Part (ii).  We have $\psi(L,P) = C$ and $u\in C\cap L$.  Let $M\le L$
be a minimal Levi subgroup containing $u$.  We wish to show that
$L=M$.  By definition, $C$ is obtained by extending to $G$ the
Richardson class in $L$ of $P$.  If $v\in L$ represents this
Richardson class in $L$ then $v$ is distinguished in $L$ by
Lemma~\ref{rep_of_dense_orbit_dist}.  Since $u$ is conjugate to $v$
(in $G$) we have $\rank \C G(u) = \rank \C G(v)$.  By
Lemma~\ref{min_Levi_unique_up_to_conj} we have $\dim Z(M) = \rank \C
G(u) = \rank \C G(v) = \dim Z(L)$ whence $L=M$.

Part (iii).  Corollary~\ref{richardson_map_inj_on_dist} shows that
$\psi$, restricted to those pairs where $L=G$, is injective and part
(ii) shows that the image of this restriction is a subset of the
distinguished classes of $G$.  Then
Corollary~\ref{bijection_reduction_to_dist_classes} shows that $\psi$
defined on all of $\BCpairs(G)$ is injective.

For surjectivity, we have two cases.  If $G$ is a classical group, we
use the description of distinguished unipotent classes in
\cite[II.7.10]{spaltenstein} or \cite{duckworth} and apply
Lemma~\ref{partitions_of_distinguished_richardson_classes} to see that
$\psi$, applied to those pairs $(L,P)$ where $L=G$, has image equal to
all the distinguished classes of $G$.  Then
Corollary~\ref{bijection_reduction_to_dist_classes} shows that $\psi$
is surjective.  If $G$ is exceptional it is simpler to count all pairs
$(L,P)$ and compare this to the number of unipotent classes in $G$ as
found in \cite{lawther2}, which draws on \cite{mizuno1},
\cite{mizuno2}, \cite{shinoda}, \cite{shoji}, \cite{G_2}.
\end{numproof}

\end{document}